\documentclass[12pt]{article}

\usepackage{amsmath, epsfig, cite}
\usepackage{amssymb}
\usepackage{amsfonts}
\usepackage{latexsym}
\usepackage{float}
\usepackage{color}

\newtheorem{thm}{Theorem}[section]

\newtheorem{lem}[thm]{Lemma}

\newcommand*{\fl}[2]{\left\lfloor\frac{#1}{#2}\right\rfloor}

\numberwithin{equation}{section}

\newcommand{\qed}{{\hfill$\square$}\medskip}

\begin{document}

\begin{center}
{\Large\bf On two congruences involving Franel numbers}
\end{center}

\vskip 2mm \centerline{Ji-Cai Liu}
\begin{center}
{\footnotesize Department of Mathematics, Wenzhou University, Wenzhou 325035, PR China\\
{\tt jcliu2016@gmail.com } \\[10pt]
}
\end{center}


\vskip 0.7cm \noindent{\bf Abstract.}
Via symbolic summation method, we establish the following series for $\pi^2$:
\begin{align*}
\sum_{k=1}^{\infty}\frac{H_k-2H_{2k}}{(-3)^k k}=\frac{\pi^2}{18},
\end{align*}
where $H_k=\sum_{j=1}^k 1/j$. We also derive a $p$-adic congruence related to this series.
As an application, we prove two congruences involving Franel numbers, one of which was originally
conjectured by Sun.


\vskip 2mm
\noindent{\it MR Subject Classifications}: 11A07, 05A19, 33F10

\section{Introduction}
It is well-known that (see \cite[(3)]{poorten-mi-1978})
\begin{align}
\sum_{k=1}^{\infty}\frac{1}{k^2{2k\choose k}}=\frac{\pi^2}{18},\label{aa-1}
\end{align}
which can be derived from the following familiar power series expansion:
\begin{align*}
2\left(\arcsin(x/2)\right)^2=\sum_{k=1}^{\infty}\frac{x^{2k}}{k^2{2k\choose k}}.
\end{align*}

Recall that the Euler numbers and Bernoulli polynomials are defined as
\begin{align*}
\frac{2}{e^x+e^{-x}}=\sum_{n=0}^{\infty}E_n\frac{x^n}{n!}
\end{align*}
and
\begin{align*}
\frac{xe^{tx}}{e^x-1}=\sum_{k=0}^{\infty}B_k(t)\frac{x^k}{k!}.
\end{align*}
In 2011, Sun \cite[(1.3)]{sunzw-scm-2011} showed that \eqref{aa-1} possesses the following interesting $p$-adic analogue:
\begin{align*}
\sum_{k=1}^{(p-1)/2}\frac{1}{k^2{2k\choose k}}\equiv (-1)^{(p-1)/2}\frac{4}{3}E_{p-3}\pmod{p},
\end{align*}
for any prime $p\ge 5$. Mattarei and Tauraso \cite[Theorem 6.1]{mt-jnt-2013} investigated
congruence properties for the polynomial $\sum_{k=1}^{p-1}k^{-2}{2k\choose k}^{-1}t^k.$

In recent decades, infinite sums on harmonic numbers related to powers of $\pi$ have been
widely studied. For instance, D. Borwein and J.M. Borwein \cite{bb-pams-1995} proved that
\begin{align*}
\sum_{k=1}^{\infty}\frac{H_k^2}{k^2}=\frac{17}{360}\pi^4,
\end{align*}
where the $k$th generalized harmonic number is given by
\begin{align*}
H_k^{(r)}=\sum_{j=1}^k\frac{1}{j^r},
\end{align*}
with the convention that $H_k=H_k^{(1)}$. It is worth mentioning that Sun \cite{sunzw-a-2011,sunzw-a-2019} conjectured and proved many infinite identities involving harmonic numbers related to powers of $\pi$. For more series for $\pi$ and related congruences as well as $q$-congruences, one can refer to \cite{ccl-am-2004,cwz-ijm-2013,guo-aam-2020,gs-rm-2019,gs-ca-2019,gz-am-2019,lp-a-2019,sunzw-scm-2011,sunzw-pams-2012,wy-a-2019,zudilin-jnt-2009}.

The first aim of the paper is to establish the following series for $\pi^2$.
\begin{thm}\label{t-1}
We have
\begin{align}
\sum_{k=1}^{\infty}\frac{H_k-2H_{2k}}{(-3)^k k}=\frac{\pi^2}{18}.\label{a-1}
\end{align}
\end{thm}

Our proof of \eqref{a-1} makes use of a finite identity, which is derived from symbolic summation
method. We also show that \eqref{a-1} has the following $p$-adic analogue.
\begin{thm}\label{t-2}
For any prime $p\ge 5$, we have
\begin{align}
\sum_{k=1}^{(p-1)/2}\frac{H_k-2H_{2k}}{(-3)^k k}\equiv \frac{1}{6}\left(\frac{p}{3}\right)B_{p-2}\left(\frac{1}{3}\right)\pmod{p},\label{a-2}
\end{align}
where $\left(\frac{\cdot}{3}\right)$ denotes the Legendre symbol.
\end{thm}

In 1894, Franel \cite{franel-1984} found that the sums of cubes of binomial coefficients:
\begin{align*}
f_n=\sum_{k=0}^n{n\choose k}^3
\end{align*}
satisfy the recurrence:
\begin{align*}
(n +1)^2f_{n+1}=(7n^2+7n+2)f_n+8n^2f_{n-1}.
\end{align*}
The numbers $f_n$ are known as Franel numbers, which also
appear in Strehl's identity (see \cite{strehl-1993}):
\begin{align}
f_n=\sum_{k=0}^n{n\choose k}^2{2k\choose n}=\sum_{k=0}^n{n\choose k}{k\choose n-k}{2k\choose k}.\label{new-1}
\end{align}
Since the appearance of the Franel numbers, some interesting congruence properties have been gradually discovered (see \cite{guo-itsf,jv-rj-2010,sunzw-aam-2013,sunzw-jnt-2013}). For instance,
Sun \cite[Theorem 1.1]{sunzw-aam-2013} proved that for any prime $p\ge 5$,
\begin{align*}
\sum_{k=0}^{p-1}(-1)^kf_k\equiv \left(\frac{p}{3}\right)\pmod{p^2}.
\end{align*}

The second aim of the paper is prove the following two congruences involving Franel numbers and
harmonic numbers.
\begin{thm}\label{t-3}
For any prime $p\ge 5$, we have
\begin{align}
&\sum_{k=0}^{p-1}(-1)^kf_kH_k^{(2)}\equiv \frac{1}{2}B_{p-2}\left(\frac{1}{3}\right)\pmod{p},\label{a-3}\\[10pt]
&\sum_{k=0}^{p-1}(-1)^kf_kH_k\equiv -2\left(\frac{p}{3}\right)q_p(3) \pmod{p},\label{a-4}
\end{align}
where $q_p(3)$ is the Fermat quotient $(3^{p-1}-1)/p$.
\end{thm}

Note that \eqref{a-3} was originally conjectured by Sun \cite[Conjecture 57]{sunzw-2019}. Our proof of \eqref{a-3} is based on Theorem \ref{t-2}. We shall prove Theorems \ref{t-1}--\ref{t-3} in Sections 2--4, respectively.

\section{Proof of Theorem \ref{t-1}}
Before proving Theorem \ref{t-1}, we require the following finite identity.
\begin{lem}
For any non-negative integer $n$, we have
\begin{align}
\sum_{k=0}^{n}(-4)^k{n\choose k}\sum_{j=1}^k\frac{1}{j^2{2j\choose j}}=
(-3)^n\sum_{k=1}^n\frac{H_k-2H_{2k}}{(-3)^k k}.\label{b-1}
\end{align}
\end{lem}
{\noindent\it Proof.}
By using symbolic summation package {\tt Sigma} due to Schneider \cite{schneider-slc-2007}, we find that both sides of \eqref{b-1} satisfy the same recurrence:
\begin{align*}
&3(n+1)(2n+3)S(n)-(2n+5)(5n+9)S(n+1)+(2n^2+17n+29)S(n+2)\\
&+(n+3)(2n+5)S(n+3)=0.
\end{align*}
It is trivial to check that both sides of \eqref{b-1} are equal for $n=0,1,2$.
\qed

{\noindent\it Proof of \eqref{a-1}.}
By \eqref{b-1}, we have
\begin{align}
\sum_{k=1}^n\frac{H_k-2H_{2k}}{(-3)^k k}
&=\left(-\frac{1}{3}\right)^n\sum_{j=1}^{n}\frac{1}{j^2{2j\choose j}}\sum_{k=j}^{n}(-4)^k{n\choose k}\notag\\[10pt]
&=\left(-\frac{1}{3}\right)^n\sum_{j=1}^{n}\frac{1}{j^2{2j\choose j}}\left((-3)^n-\sum_{k=0}^{j-1}(-4)^k{n\choose k}\right)\notag\\[10pt]
&=\sum_{j=1}^{n}\frac{1}{j^2{2j\choose j}}\left(1-\sum_{k=0}^{j-1}\frac{(-4)^k}{(-3)^n}{n\choose k}\right).\label{b-2}
\end{align}
It is obvious that for $0\le k\le j-1$,
\begin{align*}
\lim_{n\to \infty}\frac{(-4)^k}{(-3)^n}{n\choose k}=0,
\end{align*}
and so
\begin{align}
\lim_{n\to \infty}\sum_{k=0}^{j-1}\frac{(-4)^k}{(-3)^n}{n\choose k}=0.\label{b-3}
\end{align}
Letting $n\to \infty$ on both sides of \eqref{b-2} and noting \eqref{b-3}, we arrive at
\begin{align}
\sum_{k=1}^{\infty}\frac{H_k-2H_{2k}}{(-3)^k k}=
\sum_{j=1}^{\infty}\frac{1}{j^2{2j\choose j}}.\label{b-4}
\end{align}
Then the proof of \eqref{a-1} follows from \eqref{aa-1} and \eqref{b-4}.
\qed

\section{Proof of Theorem \ref{t-2}}
In order to prove Theorem \ref{t-2}, we need three preliminary results.
\begin{lem}
For any non-negative integer $n$, we have
\begin{align}
{2n\choose n}\sum_{k=1}^n\frac{1}{k^2{2k\choose k}}={2n\choose n}H_n^{(2)}
+\sum_{k=1}^n\frac{-3\left(\frac{k}{3}\right)^2+2}{k^2}{2n\choose n+k}.\label{c-1}
\end{align}
\end{lem}
{\noindent\it Proof.}
An identity due to Mattarei and Tauraso \cite[(22)]{mt-jnt-2013} says
\begin{align}
{2n\choose n}\sum_{k=1}^n\frac{t^k}{k^2{2k\choose k}}={2n\choose n}H_n^{(2)}
+\sum_{k=1}^n\frac{v_k(t-2)}{k^2}{2n\choose n+k},\label{c-2}
\end{align}
where the Lucas sequences $\{v_k(t)\}_{k\ge 0}$ are defined by
\begin{align*}
v_0(t)=2,\quad v_1(t)=t,\quad\text{and}\quad v_k(t)=tv_{k-1}(t)-v_{k-2}(t)\quad \text{for $k\ge 2$}.
\end{align*}
It is not hard to check that
\begin{align}
v_k(-1)=-3\left(\frac{k}{3}\right)^2+2.\label{c-3}
\end{align}
Letting $t=1$ in \eqref{c-2} and noting \eqref{c-3}, we reach \eqref{c-1}.
\qed

\begin{lem}
For any prime $p\ge 5$, we have
\begin{align}
\sum_{k=0}^{(p-1)/2}{2k\choose k}H_k^{(2)}\equiv -\frac{1}{6}B_{p-2}\left(\frac{1}{3}\right)\pmod{p}.\label{c-4}
\end{align}
\end{lem}
{\noindent\it Proof.}
We begin with the following polynomial congruence \cite[(37)]{mt-jnt-2013}:
\begin{align}
\sum_{k=1}^{p-1}t^{p-k}{2k\choose k}H_k^{(2)}\equiv -2t\sum_{k=1}^{p-1}\frac{u_k(2-t)}{k^2}\pmod{p},\label{c-5}
\end{align}
where the Lucas sequences $\{u_k(t)\}_{k\ge 0}$ are defined by
\begin{align*}
u_0(t)=0,\quad u_1(t)=1,\quad\text{and}\quad u_k(t)=tu_{k-1}(t)-u_{k-2}(t)\quad \text{for $k\ge 2$}.
\end{align*}
It is trivial to check that
\begin{align}
u_k(1)=\frac{(-1)^{\fl{k}{3}}+(-1)^{\fl{k-1}{3}}}{2},\label{c-6}
\end{align}
where $\lfloor x \rfloor$ denotes the integral part of real $x$.
Letting $t=1$ in \eqref{c-5} and using \eqref{c-6} gives
\begin{align}
\sum_{k=1}^{p-1}{2k\choose k}H_k^{(2)}\equiv -\sum_{k=1}^{p-1}\frac{(-1)^{\fl{k}{3}}+(-1)^{\fl{k-1}{3}}}{k^2}\pmod{p}.\label{c-7}
\end{align}

Note that
\begin{align}
&\sum_{k=1}^{p-1}\frac{(-1)^{\fl{k}{3}}+(-1)^{\fl{k-1}{3}}}{k^2}\notag\\[10pt]
&=2\left(\sum_{\substack{1\le k \le p-1\\[3pt]k\equiv1\pmod{6}}}\frac{1}{k^2}
+\sum_{\substack{1\le k \le p-1\\[3pt]k\equiv2\pmod{6}}}\frac{1}{k^2}
-\sum_{\substack{1\le k \le p-1\\[3pt]k\equiv4\pmod{6}}}\frac{1}{k^2}
-\sum_{\substack{1\le k \le p-1\\[3pt]k\equiv5\pmod{6}}}\frac{1}{k^2}\right)\notag\\[10pt]
&=2H_{p-1}^{(2)}-2\sum_{\substack{1\le k \le p-1\\[3pt]k\equiv0\pmod{3}}}\frac{1}{k^2}-4\left(\sum_{\substack{1\le k \le p-1\\[3pt]k\equiv4\pmod{6}}}\frac{1}{k^2}
+\sum_{\substack{1\le k \le p-1\\[3pt]k\equiv5\pmod{6}}}\frac{1}{k^2}\right)\notag\\[10pt]
&=2H_{p-1}^{(2)}-\frac{2}{9}H_{\lfloor p/3\rfloor}^{(2)}
-\frac{1}{9}\left(\sum_{k=1}^{\fl{p+1}{6}}\frac{1}{(k-1/3)^2}+\sum_{k=1}^{\fl{p}{6}}\frac{1}{(k-1/6)^2}\right).
\label{cc-1}
\end{align}

If $p\equiv 1\pmod{6}$, then $\lfloor p/6\rfloor\equiv -1/6\pmod{p}$, and so
\begin{align}
&\sum_{k=1}^{\fl{p+1}{6}}\frac{1}{(k-1/3)^2}+\sum_{k=1}^{\fl{p}{6}}\frac{1}{(k-1/6)^2}\notag\\[10pt]
&\equiv \sum_{k=1}^{\fl{p}{6}}\frac{1}{(2\lfloor p/6\rfloor+k)^2}+\sum_{k=1}^{\fl{p}{6}}\frac{1}{(\lfloor p/6\rfloor+k)^2}\pmod{p}\notag\\[10pt]
&=H_{(p-1)/2}^{(2)}-H_{\lfloor p/6 \rfloor}^{(2)}.\label{c-10}
\end{align}

It follows from \eqref{cc-1} and \eqref{c-10} that
\begin{align}
\sum_{k=1}^{p-1}\frac{(-1)^{\fl{k}{3}}+(-1)^{\fl{k-1}{3}}}{k^2}
\equiv 2H_{p-1}^{(2)}-\frac{2}{9}H_{\lfloor p/3\rfloor}^{(2)}
-\frac{1}{9}H_{(p-1)/2}^{(2)}+\frac{1}{9}H_{\lfloor p/6 \rfloor}^{(2)}\pmod{p}.
\label{c-13}
\end{align}
By Wolstenholme's theorem \cite[page 114]{hw-b-2008}, \cite[Lemma 2.4]{sunzw-scm-2011}, \cite[(9)]{lehmer-am-1938} and the congruence below equation \cite[(47)]{lehmer-am-1938}, we have
\begin{align}
&H_{p-1}^{(2)}\equiv 0\pmod{p},\label{cc-2}\\
&H_{(p-1)/2}^{(2)}\equiv 0\pmod{p},\label{c-11}\\
&H_{\lfloor p/3\rfloor}^{(2)}\equiv \frac{1}{2}\left(\frac{p}{3}\right)B_{p-2}\left(\frac{1}{3}\right)\pmod{p},\label{c-8}\\
&H_{\lfloor p/6 \rfloor}^{(2)}\equiv 5H_{\lfloor p/3 \rfloor}^{(2)}\overset{\eqref{c-8}}{\equiv} \frac{5}{2}\left(\frac{p}{3}\right)B_{p-2}\left(\frac{1}{3}\right)\pmod{p}.\label{c-12}
\end{align}
Substituting \eqref{cc-2}--\eqref{c-12} into \eqref{c-13} gives
\begin{align}
\sum_{k=1}^{p-1}\frac{(-1)^{\fl{k}{3}}+(-1)^{\fl{k-1}{3}}}{k^2}
\equiv \frac{1}{6}B_{p-2}\left(\frac{1}{3}\right)\pmod{p}.
\label{cc-3}
\end{align}

If $p\equiv 5\pmod{6}$, in a similar way, we can also prove that
\begin{align}
\sum_{k=1}^{p-1}\frac{(-1)^{\fl{k}{3}}+(-1)^{\fl{k-1}{3}}}{k^2}
\equiv \frac{1}{6}B_{p-2}\left(\frac{1}{3}\right)\pmod{p}.
\label{c-14}
\end{align}

Finally, combining \eqref{c-7}, \eqref{cc-3} and \eqref{c-14}, we arrive at
\begin{align*}
\sum_{k=1}^{p-1}{2k\choose k}H_k^{(2)}\equiv -\frac{1}{6}B_{p-2}\left(\frac{1}{3}\right)\pmod{p},
\end{align*}
as desired.
\qed

\begin{lem}
Let $p\ge 5$ be a prime. Then for $0\le d \le \frac{p-1}{2}$, we have
\begin{align}
\sum_{k=0}^{(p-1)/2}{2k\choose k+d}\equiv \left(\frac{p-d}{3}\right)-\left(\frac{d}{3}\right)\pmod{p}.\label{c-15}
\end{align}
\end{lem}
{\noindent\it Proof.}
From \cite[Theorem 4.2]{tauraso-aam-2012}, we deduce that for $0\le d\le n$,
\begin{align*}
\sum_{k=0}^{n-1}{2k\choose k+d}=\sum_{k=0}^{n-d}\left(\frac{n-d-k}{3}\right){2n\choose k}.
\end{align*}
Letting $n=(p+1)/2$ in the above gives
\begin{align*}
\sum_{k=0}^{(p-1)/2}{2k\choose k+d}=\sum_{k=0}^{(p+1)/2-d}\left(\frac{(p+1)/2-d-k}{3}\right){p+1\choose k}.
\end{align*}
Since for $2\le k\le (p+1)/2-d$, we have
\begin{align*}
{p+1\choose k}\equiv 0\pmod{p},
\end{align*}
and so
\begin{align}
\sum_{k=0}^{(p-1)/2}{2k\choose k+d}\equiv \left(\frac{(p+1)/2-d}{3}\right)+
\left(\frac{(p+1)/2-d-1}{3}\right)\pmod{p}.\label{c-16}
\end{align}

Next, we distinguish two cases $p\equiv 1\pmod{6}$ and $p\equiv 5\pmod{6}$ to prove \eqref{c-15}.
In both cases, we can show that
\begin{align}
\left(\frac{(p+1)/2-d}{3}\right)+
\left(\frac{(p+1)/2-d-1}{3}\right)\equiv \left(\frac{p-d}{3}\right)-
\left(\frac{d}{3}\right)\pmod{p}.\label{c-17}
\end{align}

Finally, combining \eqref{c-16} and \eqref{c-17}, we complete the proof of \eqref{c-15}.
\qed

Now we are ready to prove Theorem \ref{t-2}.

{\noindent\it Proof of \eqref{a-2}.}
Letting $n=\frac{p-1}{2}$ in \eqref{b-1} and noting that
\begin{align*}
{(p-1)/2\choose k}\equiv \frac{{2k\choose k}}{(-4)^k}\pmod{p}
\quad\text{and}\quad (-3)^{(p-1)/2}\equiv \left(\frac{p}{3}\right)\pmod{p},
\end{align*}
we obtain
\begin{align*}
\sum_{k=0}^{(p-1)/2}{2k\choose k}\sum_{j=1}^k\frac{1}{j^2{2j\choose j}}\equiv
\left(\frac{p}{3}\right)\sum_{k=1}^{(p-1)/2}\frac{H_k-2H_{2k}}{(-3)^k k}\pmod{p}.
\end{align*}
In order to prove \eqref{a-2}, it suffices to show that
\begin{align}
\sum_{k=0}^{(p-1)/2}{2k\choose k}\sum_{j=1}^k\frac{1}{j^2{2j\choose j}}
\equiv \frac{1}{6}B_{p-2}\left(\frac{1}{3}\right)\pmod{p}.\label{c-18}
\end{align}

By \eqref{c-1}, \eqref{c-4} and \eqref{c-15}, we have
\begin{align}
&\sum_{k=0}^{(p-1)/2}{2k\choose k}\sum_{j=1}^k\frac{1}{j^2{2j\choose j}}\notag\\[10pt]
&\equiv -\frac{1}{6}B_{p-2}\left(\frac{1}{3}\right)
+\sum_{k=0}^{(p-1)/2}\sum_{j=1}^k\frac{-3\left(\frac{j}{3}\right)^2+2}{j^2}{2k\choose k+j}\pmod{p}\notag\\[10pt]
&=-\frac{1}{6}B_{p-2}\left(\frac{1}{3}\right)
+\sum_{j=1}^{(p-1)/2}\frac{-3\left(\frac{j}{3}\right)^2+2}{j^2}\sum_{k=0}^{(p-1)/2}{2k\choose k+j}\notag\\[10pt]
&\equiv -\frac{1}{6}B_{p-2}\left(\frac{1}{3}\right)
+\sum_{j=1}^{(p-1)/2}\frac{\left(-3\left(\frac{j}{3}\right)^2+2\right)\left(\left(\frac{p-j}{3}\right)-\left(\frac{j}{3}\right)\right)}{j^2}\pmod{p}.
\label{c-19}
\end{align}

If $p\equiv 1\pmod{6}$, then $\lfloor p/6\rfloor\equiv -1/6\pmod{p}$, and so
\begin{align}
&\sum_{j=1}^{(p-1)/2}\frac{\left(-3\left(\frac{j}{3}\right)^2+2\right)\left(\left(\frac{p-j}{3}\right)-\left(\frac{j}{3}\right)\right)}{j^2}\notag\\[10pt]
&=\sum_{\substack{1\le j \le (p-1)/2\\[3pt]j\equiv1\pmod{3}}}\frac{1}{j^2}
+2\sum_{\substack{1\le j \le (p-1)/2\\[3pt]j\equiv0\pmod{3}}}\frac{1}{j^2}\notag\\[10pt]
&=H_{(p-1)/2}^{(2)}+\frac{1}{9}H_{\lfloor p/6\rfloor}^{(2)}-\sum_{\substack{1\le j \le (p-1)/2\\[3pt]j\equiv2\pmod{3}}}\frac{1}{j^2}.\label{c-20}
\end{align}
Note that
\begin{align}
\sum_{\substack{1\le j \le (p-1)/2\\[3pt]j\equiv2\pmod{3}}}\frac{1}{j^2}
&=\frac{1}{9}\sum_{j=1}^{\fl{p}{6}}\frac{1}{(j-1/3)^2}\notag\\
&\equiv \frac{1}{9}\sum_{j=1}^{\fl{p}{6}}\frac{1}{(j+2\lfloor p/6\rfloor)^2}\pmod{p}\notag\\[10pt]
&=\frac{1}{9}\left(H_{(p-1)/2}^{(2)}-H_{\lfloor p/3\rfloor}^{(2)}\right).\label{c-21}
\end{align}
It follows from \eqref{c-11}, \eqref{c-8}, \eqref{c-12}, \eqref{c-20} and \eqref{c-21} that
\begin{align}
\sum_{j=1}^{(p-1)/2}\frac{\left(-3\left(\frac{j}{3}\right)^2+2\right)\left(\left(\frac{p-j}{3}\right)-\left(\frac{j}{3}\right)\right)}{j^2}
&\equiv \frac{8}{9}H_{(p-1)/2}^{(2)}+\frac{1}{9}H_{\lfloor p/3\rfloor}^{(2)}+\frac{1}{9}H_{\lfloor p/6\rfloor}^{(2)}\notag\\[10pt]
&\equiv \frac{1}{3}B_{p-2}\left(\frac{1}{3}\right)\pmod{p}.\label{c-22}
\end{align}

If $p\equiv 5\pmod{6}$, in a similar way, we can also prove that
\begin{align}
\sum_{j=1}^{(p-1)/2}\frac{\left(-3\left(\frac{j}{3}\right)^2+2\right)\left(\left(\frac{p-j}{3}\right)-\left(\frac{j}{3}\right)\right)}{j^2}
\equiv \frac{1}{3}B_{p-2}\left(\frac{1}{3}\right)\pmod{p}.\label{c-23}
\end{align}

Finally, combining \eqref{c-19}, \eqref{c-22} and \eqref{c-23}, we reach \eqref{c-18}.
\qed

\section{Proof of Theorem \ref{t-3}}
To prove Theorem \ref{t-3}, we require the following two identities.
\begin{lem}
For any non-negative integer $n$, we have
\begin{align}
&\sum_{k=n}^{2n}(-1)^k{k\choose n}{n\choose k-n}H_k^{(2)}=3\sum_{k=1}^n\frac{1}{k^2{2k\choose k}},\label{d-1}\\[10pt]
&\sum_{k=n}^{2n}(-1)^k{k\choose n}{n\choose k-n}H_k=2H_n.\label{d-2}
\end{align}
\end{lem}
{\noindent\it Proof.}
The identities \eqref{d-1} and \eqref{d-2} possess the same proofs as \eqref{b-1}, and we omit the details.
\qed

{\noindent\it Proof of \eqref{a-3}.}
By \eqref{new-1}, we have
\begin{align}
\sum_{k=0}^{p-1}(-1)^kf_kH_k^{(2)}
&=\sum_{k=0}^{p-1}(-1)^kH_k^{(2)}\sum_{j=0}^k{k\choose j}{j\choose k-j}{2j\choose j}\notag\\
&=\sum_{j=0}^{p-1}{2j\choose j}\sum_{k=0}^{p-1}(-1)^k{k\choose j}{j\choose k-j}H_{k}^{(2)}\notag\\
&\equiv \sum_{j=0}^{(p-1)/2}{2j\choose j}\sum_{k=0}^{p-1}(-1)^k{k\choose j}{j\choose k-j}H_{k}^{(2)}\pmod{p},\label{d-3}
\end{align}
because ${2j\choose j}\equiv 0\pmod{p}$ for $(p+1)/2\le j\le p-1$.

Note that for $0\le j\le (p-1)/2$,
\begin{align}
\sum_{k=0}^{p-1}(-1)^k{k\choose j}{j\choose k-j}H_{k}^{(2)}&=\sum_{k=0}^{2j}(-1)^k{k\choose j}{j\choose k-j}H_{k}^{(2)}\notag\\
&\overset{\eqref{d-1}}{=}3\sum_{k=1}^j\frac{1}{k^2{2k\choose k}}.\label{d-4}
\end{align}
It follows from \eqref{d-3}, \eqref{d-4} and \eqref{c-18} that
\begin{align*}
\sum_{k=0}^{p-1}(-1)^kf_kH_k^{(2)}&\equiv
3\sum_{j=0}^{(p-1)/2}{2j\choose j}\sum_{k=1}^j\frac{1}{k^2{2k\choose k}}\\
&\equiv \frac{1}{2}B_{p-2}\left(\frac{1}{3}\right)\pmod{p},
\end{align*}
as desired.
\qed

{\noindent\it Proof of \eqref{a-4}.}
Similarly to the proof of \eqref{a-3}, by using \eqref{d-2}, we can show that
\begin{align*}
\sum_{k=0}^{p-1}(-1)^kf_kH_k\equiv 2\sum_{j=0}^{(p-1)/2}{2j\choose j}H_j\pmod{p}.
\end{align*}
Combining the above and the following congruence (see \cite[page 528]{ms-ijnt-2016}):
\begin{align*}
\sum_{j=0}^{(p-1)/2}{2j\choose j}H_j\equiv -\left(\frac{p}{3}\right)q_p(3)\pmod{p},
\end{align*}
we complete the proof of \eqref{a-4}.
\qed

\vskip 5mm \noindent{\bf Acknowledgments.}
This work was supported by the National Natural Science Foundation of China (grant 11801417).

\end{document}